\documentclass[reqno]{amsart} \usepackage{amscd}
\usepackage{epsf}
\newtheorem{theorem}{Theorem}[section]
\newtheorem{proposition}[theorem]{Proposition}

\theoremstyle{remark} \newtheorem{remark}[theorem]{Remark}

\usepackage{graphicx}

\begin{document}

	\title[On Euler-Ramanujan formula]{On Euler-Ramanujan formula, Dirichlet series and minimal surfaces}
	
	\author{Rukmini Dey$^1$, Rishabh Sarma$^2$, Rahul Kumar Singh$^3$}
	\address{1.  I.C.T.S.-T.I.F.R.\\
		Bangalore, India\\
		rukmini@icts.res.in\\
		rukmini.dey@gmail.com\\
		2. I.C.T.S.-T.I.F.R.\\
		Bangalore, India\\
		3. School of Mathematics, National Institute of Science Education and Research, Bhubaneswar (HBNI), Odisha 752050, India\\
		rhlsngh498@gmail.com}

	\begin{abstract}
		In this paper, we rewrite two forms of an Euler-Ramanujan identity in terms of certain Dirichlet series and derive functional equation of the latter. We also use the Weierstrass-Enneper representation of minimal surfaces to obtain some identities involving these Dirichlet series and one complex parameter. 
	\end{abstract}
	
	\maketitle
	
	\section{Introduction}
	
	In this article, we explore two different forms of the Euler-Ramanujan identities and their relationship to certain Dirichlet series and minimal surfaces.  The use of Weierstrass-Enneper representation of minimal surfaces in E-R identities first appeared  in ~\cite{dey}.  When we rewrite the identities in terms of sums of  Dirichlet series we  get relationships between the latter  and sometimes using the W-E representation of minimal surfaces the relationship is in terms of one complex parameter. The minimal surfaces which are involved are a family of Scherk's surfaces and  helicoids which are  minimal surfaces of translation. In section $2$ of this paper we devote ourselves to   relationships between sums of Dirichlet series. In section $3$ we study the properties of one  of the Dirichlet series and  derive a functional equation of the latter.

	\section{Rewriting a form of the  Euler-Ramanujan identity in terms of Dirichlet series }

	We have Euler-Ramanujan's identity, ~\cite{ramanujan}, Example $(1)$ page $38$, where 
	$X$, $A$ are complex variables, $A$ is not an odd multiple of $\pi /2$:
	\begin{eqnarray*}
		\frac{\cos{(X+A)}}{\cos{A}} = \prod_{k=1}^{\infty} \left( 1-\frac{X}{(k - \frac{1}{2})\pi -A} \right) \left( 1 + \frac{X}{(k - \frac{1}{2})\pi +A}\right).
	\end{eqnarray*}
	We take complex logarithm on both sides, to get:

	\begin{eqnarray*}
		\log \left(\frac{\cos{(X+A)}}{\cos{A}}\right)=\sum_{k=1}^{\infty} \log \left\lbrace \left(\frac{ (k - \frac{1}{2})\pi -(X + A)}{(k - \frac{1}{2}) \pi - A} \right) \left(\frac{ (k - \frac{1}{2}) \pi + (X + A) }{(k - \frac{1}{2}) \pi +A} \right)\right\rbrace.
	\end{eqnarray*}
	
	Let $c_k = (k-\frac{1}{2}) \pi$. Replace $X+A$ by $y$ and $A$ by $x$ in the identity  where $x$ (not an odd multiple of $\frac{\pi}{2}$) and $y$ are real parameters, we have, 
	\begin{eqnarray}\label{lnyx}
	\log \left(\frac{\cos{y}}{\cos{x}}\right) &=& \sum_{k=1}^{\infty} \log \left(\frac{c_k - y}{c_k - x} \right) \left(\frac{c_k + y}{c_k + x} \right)
	\end{eqnarray}
	
	or,
	$$\log \left(\frac{\cos{y}}{\cos{x}}\right)=\sum_{k=1}^{\infty}  [\log ( 1 + \frac{x^2}{c_k^2 - x^2} ) + \log ( 1 - \frac{y^2}{c_k^2}) ].$$
	Let us take the special case when $|\frac{x^2}{c_k^2 - x^2}| < 1$ and $|\frac{y^2}{c_k^2}| < 1$.
	
	Then, 
	\begin{eqnarray*}
		\log \left(\frac{\cos{y}}{\cos{x}}\right)
		&=& \sum_{k=1}^{\infty}  [\log ( 1 + \frac{x^2}{c_k^2 - x^2} ) + \log ( 1 - \frac{y^2}{c_k^2}) ] \\
		&=& \sum_{k=1}^{\infty} [ (\frac{x^2}{c_k^2 - x^2} - \frac{x^4}{2(c_k^2 - x^2)^2} + \frac{x^6}{3(c_k^2 - x^2)^3} -...) + ( -\frac{y^2}{c_k^2} - \frac{y^4}{2c_k^4} - \frac{y^6}{3c_k^6}...)]\\
		&=& \sum_{k=1}^{\infty} [ L_k (1, \frac{x^2}{c_k^2 -x^2}) - M_k (1, \frac{y^2}{c_k^2}) ]
	\end{eqnarray*}
	where $L_k (s, a) = \sum_n ((-1)^{n+1}(\frac{x^2}{c_k^2 -x^2})^n / n^s)$ is a Dirichlet series with a real paramater $a = \frac{x^2}{c_k^2 -x^2}$ and an integer $k$, $ M_k (s, b) = \sum_n ( \frac{y^2}{c_k^2})^n / n^s$ is a Dirichlet series with a real parameter  $b = \frac{y^2}{c_k^2}$ and an integer $k$ and  $L_k (1, \frac{x^2}{c_k^2 -x^2})$ and $ M_k (1, \frac{y^2}{c_k^2})$ are these Dirichlet series  evaluated at $s=1$.
	We note that $|\frac{x^2}{c_k^2 - x^2}| < 1$ and $|\frac{y^2}{c_k^2}| < 1$ for all $k$ is equivalent to 
	$ -\frac{\pi}{2 \sqrt{2}}  < x < \frac{\pi}{2\sqrt{2}} $ and $-\sqrt{\pi/2} < y < \sqrt{\pi/2}$.

	Thus we have the following proposition:
	\begin{proposition}
		$\log \left(\frac{\cos{y}}{\cos{x}}\right) = \sum_{k=1}^{\infty} [ L_k (1, \frac{x^2}{c_k^2 -x^2}) -M_k (1, \frac{y^2}{c_k^2}) ]$  for $-\frac{\pi}{2 \sqrt{2}}  < x < \frac{\pi}{2\sqrt{2}}$ and $ -\sqrt{\pi/2}< y < \sqrt{\pi/2}$.
	\end{proposition}
	
	In ~\cite{dey}, using the Weierstrass-Enneper representation of minimal surfaces we derived the following way of writing the equation $z = \log \left(\frac{\cos{y}}{\cos{x}}\right)$ (Scherk's minimal surface) in parametric form (in terms of a complex parameter $\zeta$),
	$x(\zeta, \bar{\zeta} ) =  2 {\rm Re} \tan^{-1}(\zeta),\; y(\zeta, \bar{\zeta}) = - {\rm Im} \log \left(\frac{1 + \zeta}{1-\zeta}\right)$ and 
	$z(\zeta, \bar{\zeta}) =  {\rm Re} \log \left(\frac {1 + \zeta^2}{ 1 - \zeta^2}\right) $.
	This parametrization fails precisely  at $\zeta = \pm 1, \pm i$.
	
	Using the log version of the E-R identity , for $\zeta \neq \pm 1, \pm i$ and belonging to a small domain in ${\mathbb C}$,  the expression in terms of Dirichlet series look as follows:

	\begin{eqnarray*}
		& & {\rm Re}\log \left(\frac{ 1 + \zeta^2}{  1 - \zeta^2}\right)\\
		&=& \sum_{k=1}^{\infty}  [\log ( 1 + \frac{(2 {\rm Re} \tan^{-1}(\zeta))^2}{c_k^2 - (2 {\rm Re} \tan^{-1}(\zeta))^2} ) + \log ( 1 - \frac{(- {\rm Im} \log \left(\frac{1 + \zeta}{1-\zeta}\right))^2}{c_k^2}) ]\\
		&=& \sum_{k=1}^{\infty} [ L_k (1, \frac{(2 {\rm Re} \tan^{-1}(\zeta))^2}{c_k^2 - (2 {\rm Re} \tan^{-1}(\zeta))^2}) - M_k (1, \frac{(- {\rm Im} \log \left(\frac{1 + \zeta}{1-\zeta}\right))^2}{c_k^2}) ].
	\end{eqnarray*}
	
	The condition that $|\frac{x^2}{c_k^2 - x^2}| < 1$ and $|\frac{y^2}{c_k^2}| < 1$ for all $k$ is satisfied if $|\zeta| < \frac{1}{2}$.
	This can be seen as follows. For $\zeta \in \mathbb{C}$ such that $|\zeta| < \frac{1}{2}$,  we have
	$ -2 < 2 \frac{|\zeta| \cos {\theta}}{1 - |\zeta|^2} < 2,$ which implies 
	$ -\tan{(\frac{\pi}{2 \sqrt{2}})} <  2 \frac{|\zeta| \cos { \theta}}{1 - |\zeta|^2} < \tan{(\frac{\pi}{2 \sqrt{2}})}. $ 
	In other words, $ -\tan{(\frac{\pi}{2 \sqrt{2}})} <  \frac{\zeta + \bar{\zeta}}{1 - |\zeta|^2} < \tan{(\frac{\pi}{2 \sqrt{2}})} $ and hence 
	$-\frac{\pi}{2 \sqrt{2}} < \tan^{-1}(\zeta) + \tan^{-1}(\bar{\zeta}) < \frac{\pi}{2 \sqrt{2}},$ or  $-\frac{\pi}{2\sqrt{2}} < 2 \rm{Re} \tan^{-1}(\zeta) < \frac{\pi}{2 \sqrt{2}}.$
	This implies  $|\frac{x^2}{c_k^2 - x^2}| < 1$ where $x = 2 \rm{Re} \tan^{-1}(\zeta)$. From this it follows that $|\frac{x^2}{c_k^2 - x^2}| < 1$ for all $k$, $y^2/c_k^2 <1$ does not give any additional constraint on $\zeta$.
	
	Thus, we have the following proposition:
	\begin{proposition}\label{W-E-Lfns}
		For $\zeta \in {\mathbb C}$ such that $|\zeta| < \frac{1}{2}$, 
		\begin{eqnarray*} 
			{\rm Re} \log \left(\frac{ 1 + \zeta^2}{  1 - \zeta^2}\right)&=& \sum_{k=1}^{\infty} [ L_k (1, \frac{(2 {\rm Re} \tan^{-1}(\zeta))^2}{c_k^2 - (2 {\rm Re} \tan^{-1}(\zeta))^2}) -M_k (1, \frac{(- {\rm Im} \log \left(\frac{1 + \zeta}{1-\zeta}\right))^2}{c_k^2}) ].
		\end{eqnarray*}
	\end{proposition}
	
	\subsection{More examples of minimal  surfaces and their relation to Dirichlet series}

	In this section, we will revisit the \textit{Scherk's} surface, \cite{scherk}, and the \textit{helicoid}, \cite{nitsche, osserman}. In fact, we consider a one parameter family of \textit{Scherk's} type surface. We show that all these surfaces will give reise to  a  Dirichlet series  decomposition. We consider the case of \textit{helicoid} separately.

	We consider a one parameter family of minimal translation surfaces (\textit{Scherk's} type surfaces),
	\begin{eqnarray*}
		X_{\theta}(u,v)=(u+v\cos{\theta},v\sin{\theta},\log{\frac{\cos v}{\cos u}}),
	\end{eqnarray*}
	where $\theta \in \mathbb{R}$. Indeed, for a fixed $\theta$, $X_{\theta}(u,v)=\alpha(u)+\beta_{\theta}(v)$, where $\alpha(u)=(u,0,-\log{\cos{u}})$ and $\beta_{\theta}(v)=(v\cos{\theta},v\sin{\theta},\log{\cos{v}})$ and hence it is a minimal surface of translation (see remark at the end of this section).
	\begin{remark}
		Observe that $\theta=0$ corresponds to a $\textit{plane}$ in the above   family  and $\theta=\frac{\pi}{2}$ corresponds to the classical $\textit{Scherk's}$ surface, i.e., $X_{\frac{\pi}{2}}(u,v)=(u,v,\log{\frac{\cos v}{\cos u}})$ (discussed in the previous section).
	\end{remark}
	
	The $\theta$-family of minimal translation surfaces gives us a family of Euler-Ramanujan's identities. For $\theta\neq \pm\frac{(2n+1)\pi}{2}$, the corresponding identity we call it a $\textit{twisted}$ Euler-Ramanujan's identity. Indeed, put $u+v\cos {\theta}=x$ and $v\sin{\theta}=y$, then using the Euler-Ramanujan's identity we obtain
	\begin{eqnarray*}
		\frac{\cos({\frac{y}{\sin{\theta}}})}{\cos(x-y\cot{\theta})}=\prod_{k=1}^{\infty}\left( \frac{c_k^2-\frac{y^2}{\sin^2{\theta}}}{c_k^2-(x-y\cot{\theta})^2}\right),
	\end{eqnarray*}
	where $x-y\cot{\theta}$ is not an odd multiple of $\frac{\pi}{2}$.
	Taking $\log$ on both the sides, we get
	\begin{eqnarray*}
		\log\left({\frac{\cos({\frac{y}{\sin{\theta}}})}{\cos(x-y\cot{\theta})}}\right)=\sum_{k=1}^{\infty}\log\left( \frac{c_k^2-\frac{y^2}{\sin^2{\theta}}}{c_k^2-(x-y\cot{\theta})^2}\right).
	\end{eqnarray*}
	When we take $\theta$ as an odd multiple of $\frac{\pi}{2}$ in the above identity, we get back the identity \eqref{lnyx}. Following the same idea (as in Scherk's surface case), we have a Dirichlet series decomposition in this case as well, which is as follows:
	\begin{proposition}
		\begin{eqnarray*}
			\log\left({\frac{\cos({\frac{y}{\sin{\theta}}})}{\cos(x-y\cot{\theta})}}\right)=\sum_{k=1}^{\infty}[ L_k(1, \frac{(x-y\cot \theta)^2}{c_k^2 -(x-y\cot \theta)^2}) - M_k(1, (\frac{y}{c_k\sin \theta })^2) ]   
		\end{eqnarray*}
		where $L_k(1, \frac{(x-y\cot \theta)^2}{c_k^2 -(x-y\cot \theta)^2})$ is a Dirichlet series (evaluated at $s=1$) with a real paramater $a= \frac{(x-y\cot \theta)^2}{c_k^2 -(x-y\cot \theta)^2}$ for a fixed $\theta$ and an integer $k$, $ M_k (s, -(\frac{y}{c_k \sin \theta })^2 )$, ($\theta$-fixed), is a Dirichlet series  (evaluated at $s=1$ ) with a real parameter $b= -(\frac{y}{c_k \sin \theta })^2$ and an integer $k$ as before.  
	\end{proposition}
	
	\vspace{.1cm}
	
	Now, let us consider the non parametric representation of helicoid which is given by $z=\tan^{-1}\frac{y}{x}$ and also recall the identity
	\begin{equation*}
	\tan^{-1}\omega=\frac{i}{2}(\log{(1-i\omega)}-\log{(1+i\omega)})
	\end{equation*}
	where $\omega$ is a complex number. 
	Now if we put $\omega=\frac{y}{x}$ in the above identity, we get
	\begin{equation}\label{tanl}
	\tan^{-1}\frac{y}{x}=\frac{i}{2}\left\lbrace\log{\left(1-i\frac{y}{x}\right)}-\log{\left(1+i\frac{y}{x}\right)}\right\rbrace  
	\end{equation}
	The above expression  helps us to write helicoid as a sum of two Dirichlet series evaluated at a specific value. When $|y|<|x|$ we can express the identity \eqref{tanl} as 
	
	\begin{eqnarray}\label{heliL}
	\tan^{-1}\frac{y}{x}&=&\frac{i}{2}\left\lbrace\sum_{k=1}^\infty                           (-1)^{2k-1}\frac{i^k}{k}\left(\frac{y}{x}\right)^k-\sum_{k=1}^\infty                       (-1)^{k+1}\frac{i^k}{k}\left(\frac{y}{x}\right)^k\right\rbrace \nonumber \\
	&=&\frac{1}{2}\left\lbrace-\sum_{k=1}^\infty \frac{i^{k+1}}{k}\left(\frac{y}{x}\right)^k+\sum_{k=1}^\infty (-1)^{k}\frac{i^{k+1}}{k}\left(\frac{y}{x}\right)^k\right\rbrace \nonumber \\
	&=&\frac{-1}{2}[L(1,\frac{y}{x}) - M(1,-\frac{y}{x})]
	\end{eqnarray}
	where $L(s, \frac{y}{x}) = \sum_{k=1}^\infty \frac{i^{k+1}}{k^s}( \frac{y}{x})^k$ and $ M (s, -\frac{y}{x}) = \sum_{k=1}^\infty \frac{i^{k+1}}{k^s} (-\frac{y}{x})^k$ are Dirichlet series  which are in turn evaluated at $s=1$.
	Thus we have the following proposition:
	
\vspace{.1cm}	

	\begin{proposition}\label{heli-L}
		When $|y|<|x|$, $\tan^{-1}\frac{y}{x}=\frac{1}{2}[L(1,\frac{y}{x})-M(1,-\frac{y}{x})]$.
	\end{proposition}

\vspace{.1cm}

	The Weierstrass-Enneper representation of the helicoid in terms of the complex parameter is given by (see for instance \cite{dey}),
	$x(\zeta, \bar{\zeta})=-\frac{1}{2}{\rm Im} \left(\zeta+\frac{1}{\zeta}\right)$, $y(\zeta, \bar{\zeta})=\frac{1}{2}{\rm Re} \left(\zeta-\frac{1}{\zeta}\right),$ and $z(\zeta, \bar{\zeta})=-\frac{\pi}{2}+{\rm Im} (\log{\zeta})$. The condition $|y| < |x|$ translates to $|\rm{Re} (\zeta - \frac{1}{\zeta}) | < | \rm{Im} (\zeta + \frac{1}{\zeta})|$ which is satisfied if $|\zeta| <1$. Now it is easy to see that we have 
	
	\begin{proposition}\label{heli-W-E-Lfns}
		For $\zeta \in {\mathbb C}$, such that $|\zeta| <1$
		\begin{equation*}
		-\frac{\pi}{2}+{\rm Im} (\log{\zeta})=\frac{-1}{2}\left\lbrace L\left(1, \frac{{\rm Re} \left(\zeta-\frac{1}{\zeta}\right)}{-{\rm Im} \left(\zeta+\frac{1}{\zeta}\right)}\right)-M\left(1,-\frac{{\rm Re} \left(\zeta-\frac{1}{\zeta}\right)}{-{\rm Im} \left(\zeta+\frac{1}{\zeta}\right)}\right)\right\rbrace.
		\end{equation*}
	\end{proposition}
	
	\vspace{.2cm}
	
	\subsection{Another Euler-Ramanujan identity}

	Next we look at another identity.
	For $ X $ and $ A $ real, we have (see entry $11$ in \cite{ramanujan})
	
	\begin{equation}\label{entry11}
	\tan^{-1}(\tanh X\cot A)=\tan^{-1}\left(\frac{X}{A}\right)+\sum_{k=1}^{\infty}\left(\tan^{-1}\left(\frac{X}{k\pi+A}\right)-\tan^{-1}\left(\frac{X}{k\pi+A}\right)\right).
	\end{equation}
	
	When $A=\frac{\pi}{2}$ the identity \eqref{entry11} reduces to the following identity
	\begin{equation}\label{heli}
	\tan^{-1}\left(\frac{2X}{\pi}\right)=\sum_{k=1}^{\infty}\left(\tan^{-1}\left(\frac{X}{c_k}\right)-\tan^{-1}\left(\frac{X}{d_k}\right)\right),  
	\end{equation}
	where  $c_k = (k-\frac{1}{2}) \pi$ and  $d_k = (k+\frac{1}{2}) \pi$. Next we put $\frac{2X}{\pi}=\frac{y}{x}$ in \eqref{heli} to obtain 
	\begin{equation}
	\tan^{-1}\left(\frac{y}{x}\right)=\sum_{k=1}^{\infty}\left(\tan^{-1}\left(\frac{e_ky}{x}\right)-\tan^{-1}\left(\frac{f_ky}{x}\right)\right),  
	\end{equation}
	where $e_k=\frac{\pi}{2c_k}$ and $f_k=\frac{\pi}{2d_k}$.
	Thus we have
	\vspace{.1cm}
	\begin{proposition}
		$$L(1,\frac{y}{x})-M(1,-\frac{y}{x}) = \sum_{k=1}^{\infty} (L_k(1, \frac{e_ky}{x}) - M_k(1, -\frac{e_ky}{x}) )-  (L_k(1, \frac{f_ky}{x} )- M_k(1, -\frac{f_ky}{x}))$$ 
		where $L_k(1, \frac{e_ky}{x} ) =L(1, \frac{e_ky}{x})$ and $M_k(1, -\frac{e_ky}{x} ) =M(1,- \frac{f_ky}{x})$.
	\end{proposition}
	
\vspace{.1cm}

	\begin{remark}
		Our objects of interest here happen to be minimal surfaces of translation. In the past, there has been a   good amount of interest in knowing what are all the minimal translation surfaces in $\mathbb{R}^3$. See for instance
		Lopez and Hasanis, ~\cite{lopez1}.
		It would be interesting to see if number-theoretical identities like the E-R identities are available for all minimal surfaces of translation.
	\end{remark}

	\section{Properties of the Dirichlet series }
	
	We consider the general form of the Dirichlet series defined in the previous section
	\begin{align*}
	&L_k(s,a)=\sum_{n=1}^{\infty} (-1)^{n-1}\frac{a^n}{n^s},~where~a=\frac{x^2}{c_k^2 -x^2}, 0<a<1,~s\in \mathbb{C}\\
	&M_k(s,a)=\sum_{n=1}^{\infty} (-1)^{n}\frac{a^n}{n^s},~where~a=\frac{y^2}{c_k^2}, 0 <a <1,~s\in \mathbb{C}.
	\end{align*}
	
	\subsection{Convergence of the series}
	
	We show that the series is convergent for all $s\in \mathbb{C}$ by showing its absolute convergence.
	\begin{align*}
	L_k(s,a)=\sum_{n=1}^{\infty} (-1)^{n-1}\frac{a^n}{n^s},0<a<1,~s=\sigma+ib,~\sigma,b \in \mathbb{R}.
	\end{align*}
	
	Now, $|n^s|=|n^\sigma||n^{ib}|=n^\sigma|e^{ib\ln n}|=n^\sigma$.
	\vspace{3mm}
	The absolute series then is 
	\begin{align*}
	\sum_{n=1}^{\infty} \frac{|a|^n}{n^\sigma} = \sum_{n=1}^{\infty} \frac{1}{|c|^nn^\sigma},~a=1/c, c>1.
	\end{align*}
	
	Then the ratio test gives 
	\begin{align*}
	\lim \limits_{n \rightarrow \infty} |\frac{a_{n+1}}{a_n}|=\lim \limits_{n \rightarrow \infty} \frac{|c|^n n^\sigma}{|c|^{n+1}(n+1)^\sigma}
	=\lim \limits_{n \rightarrow \infty} \frac{1}{|c|}\dfrac{1}{(1+1/n)^\sigma}=\frac{1}{c}<1.
	\\
	\end{align*}
	Thus $L_k(s,a)$ is convergent on $0<a<1$, for each $k$.	Similar is the case for $M_k(s,a)$.
	\subsection{Functional equation}	
	
	In this section, we follow a technique  of Hardy, expounded in $[1]$.
	
	Our general Dirichlet series is $\gamma(s)=\sum_{n=1}^{\infty} (-1)^{n-1}\frac{a^n}{n^s}$, where, $0<a<1,~s \in \mathbb{C}$.
	\\
	\\
	Then, defining $F(x)=\sum_{n \leq x}(-1)^{n-1}=\dfrac{1-{(-1)}^m}{2},~m<x<m+1$, we have,
	
	\begin{align*}
	\gamma(s)&=\sum_{n=1}^{\infty} (-1)^{n-1}\frac{a^n}{n^s}
	\\
	&=\sum_{n=1}^{\infty} [F(n)-F(n-1)]\frac{a^n}{n^s}
	\\
	&=\sum_{n=1}^{\infty} F(n)\frac{a^n}{n^s}-\sum_{n=1}^{\infty} F(n)\frac{a^{n+1}}{(n+1)^s}-F(0).a
	\\
	&=\sum_{n=1}^{\infty} F(n) \int_{n}^{n+1} \Big(\frac{a^x \ln a}{x^s}-\frac{sa^x}{x^{s+1}}\Big) dx
	\\
	&=\ln a\sum_{n=1}^{\infty} \int_{n}^{n+1} \frac{a^x F(x)}{x^s}-s\sum_{n=1}^{\infty} \int_{n}^{n+1} \frac{a^x F(x)}{x^{s+1}} dx
	\\
	&=\ln a \int_{1}^{\infty} \frac{a^x F(x)}{x^s}-s \int_{1}^{\infty} \frac{a^x F(x)}{x^{s+1}} dx.
	\end{align*}
	
	\begin{align*}
	Let~f(x)&=\ln a \int_{1}^{\infty} \frac{a^x F(x)}{x^s}
	\\
	&=\ln a \int_{1}^{\infty} \frac{a^x (F(x)-1/2)}{x^{s}} dx+\frac{\ln a}{2} \int_{1}^{\infty} \frac{a^x}{x^{s}} dx
	\\
	&=I_1 (s)+ I_2(s) ~(say),~and
	\\
	\\
	g(x)&=s \int_{1}^{\infty} \frac{a^x F(x)}{x^{s+1}} dx
	\\
	&=s \int_{1}^{\infty} \frac{a^x (F(x)-1/2)}{x^{s+1}} dx+\frac{s}{2} \int_{1}^{\infty} \frac{a^x}{x^{s+1}} dx
	\\
	&=I_3(s)  + I_4(s) ~(say).
	\end{align*}
	The integrals can all be seen to be convergent for all $s$. 
	\\
	\\
	$F(x)-\frac{1}{2}$ is bounded, say by $M$. 
	Then, taking $a=\frac{1}{c}, c>1$ we can see that $\int_{1}^{\infty} |\frac{a^x (F(x)-\frac{1}{2})}{x^{s}}|dx = \int_{1}^{\infty} |\frac{(F(x)-\frac{1}{2})}{|c|^xx^{\sigma}}|dx \leq M \int_{1}^{\infty} \frac{1}{c^xx^{\sigma}} dx$ (for $\sigma > 0$), which is convergent. $I_1$ and similarly $I_2$ are thus absolutely convergent since $c^x$ is exponential growth and $x^{\sigma}$ has only polynomial growth.  $I_3$ and $I_4$ can similarly be concluded to be convergent.
	\\
	\\
	We try to evaluate $I_4$ first. Let $x\ln a = u \Rightarrow du=\ln a~dx$. Then $u$ ranges from $\ln a$ to $-\infty$ as $x$ ranges from $1$ to $\infty$.
	
	Then,
	\begin{align*}
	I_4(s) &=\frac{s}{2} \int_{1}^{\infty} \frac{a^x}{x^{s+1}} dx
	\\
	&=\frac{s}{2} \int_{1}^{\infty} x^{-s-1}e^{x\ln a} dx
	\\
	&=\frac{s}{2\ln a} \int_{\ln a}^{-\infty} (\frac{u}{\ln a})^{-s-1} e^u du
	\\
	&=\frac{s(\ln a)^s}{2} \int_{\ln a}^{-\infty} u^{-s-1} e^u du
	\\
	&=\frac{-s(\ln a)^s}{2} \int_{-\ln a}^{\infty} (-1)^{-s-1}u^{-s-1} e^{-u} du
	\\
	&=\frac{s}{2} (-\ln a)^s~\Gamma(-s,-\ln a)
	\end{align*}
	
	Similar to $I_4(s)$, $I_2(s) $ will then be $-\frac{(-\ln(a))^s }{2}~\Gamma (-(s-1),-\ln a)$, where $\Gamma$ is the incomplete gamma function.
	\\
	\\
	We now calculate $I_3$. Since the function $F(x)$ is piecewise continuous of period 2, we calculate its Fourier expansion. Then, $F(x)=\frac{a_0}{2} + \sum_{n=1}^{\infty} (a_ncos\frac{n\pi x}{L}+b_nsin\frac{{n\pi x}}{L})$, where, 
	
	$a_0=\frac{1}{L}\int_{0}^{2L} F(x)cos\frac{0\pi x}{L}dx =\int_{0}^{2} F(x)dx=\int_{1}^{2} dx=1$

	$a_m==\frac{1}{L}\int_{0}^{2L} F(x)cos\frac{m\pi x}{L}dx =0$
	
	$b_n=\frac{1}{L}\int_{0}^{2L} F(x)sin\frac{n\pi x}{L} dx= \frac{-2}{n\pi}$ for $n$ odd.
	
	Therefore, $F(x)=\frac{1}{2}+\sum_{k=0}^{\infty} \dfrac{-2sin(2k\pi x+\pi x)}{(2k+1)\pi}$.\\\\\\
	Substituting $F(x)-\frac{1}{2}$ by its Fourier expansion, we have,
	\begin{align*}
	I_3(s)&=s \int_{1}^{\infty} \frac{a^x (F(x)-1/2)}{x^{s+1}} dx\\
	&=\dfrac{-2s}{\pi} \sum_{k=0}^{\infty}\dfrac{1}{2k+1}\int_{0}^{\infty} \frac{a^x sin(2k\pi x+\pi x)}{x^{s+1}} dx\\
	&=\dfrac{-2s}{\pi} \sum_{k=0}^{\infty}\dfrac{1}{2k+1}\int_{0}^{\infty} \frac{a^{\frac{t}{\pi (2k+1)}} \sin(t)}{\Big(\dfrac{t}{\pi (2k+1)}\Big)^{s+1}}~\dfrac{dt}{(2k+1)\pi}~(Taking~2k \pi x+ \pi x=t)\\
	&=\dfrac{-2s}{\pi} \sum_{k=0}^{\infty}\dfrac{\pi^s}{(2k+1)^{1-s}}\int_{0}^{\infty} \frac{a^{\frac{t}{\pi (2k+1)}} \sin(t)}{t^{s+1}}~dt.\
	\end{align*}
	We now try to evaluate $I_{5k}(s)=\int_{0}^{\infty} \frac{a^{\frac{t}{\pi (2k+1)}} \sin(t)}{t^{s+1}}~dt$.

	Let $\rm{Re}(s) <0$.

	$I_{5k}$ can be evaluated in terms of a complex gamma function as follows:
	
	\begin{eqnarray*}
		I_{5k} (s)&=& \int_{0}^{\infty} t^{-s-1} a^{\frac{t}{\pi(2k+1)}} sin(t) =  \int_{0}^{\infty} t^{-s-1} a^{\frac{t}{\pi(2k+1)}} (e^{it} - e^{-it})/2i \\
		&=& \int_{0}^{\infty} t^{-s-1} (e^{C_+ t} - e^{C_-t})/2i 
	\end{eqnarray*}
	where $C_{\pm} = \frac{ln(a)}{ \pi (2k+1)} \pm  i $.

	Let $u = C_+ t $ and $w=C_- t$. Then making change to these complex variables  and multiplying by $(-1)^{-s -1}$ when needed, we get
	
	\begin{eqnarray*}
		I_{5k}(s) = \frac{1}{2i C_{+}}  \int_{\gamma_1} u^{-s-1} e^{-u} du - \frac{1}{2 i C_{-}} \int_{\gamma_2} w^{-s -1} e^{-w} dw 
	\end{eqnarray*}
	where $\gamma_1$ is  from $0$ to $C_{+} \infty$ along the line $y = \frac{(2k+1) \pi}{ln(a)} x$ and  $\gamma_2$ is from $0$ to $C_{-}  \infty$ along the  line  $y = \frac{-(2k+1) \pi}{ln(a)} x$. Take contours shown in the diagrams below (with $\gamma_1 =e$ in the first diagram and $\gamma_2= e$ in the second diagram) and let the radius of the circular arc grow bigger.  It can be shown that the first  integral evaluates to $ \int_{\infty}^0  u^{-s -1} e^{-u} du = -\Gamma(-s) $ and so does the second one, if ${\rm Re}(s)<0$. Here  $\Gamma(-s) \equiv \int_{0}^{\infty} u^{-s -1} e^{-u} du $ and has poles along non-positive integers for $\rm{Re}(s) <0$.
	It is easy to check that $I_{5k}(s)= -\Gamma(-s)/ ((\frac{\ln(a)}{(2k+1)\pi})^2 + 1)$.

	\begin{figure}[htb!]
		\centering{\includegraphics[scale=0.6]{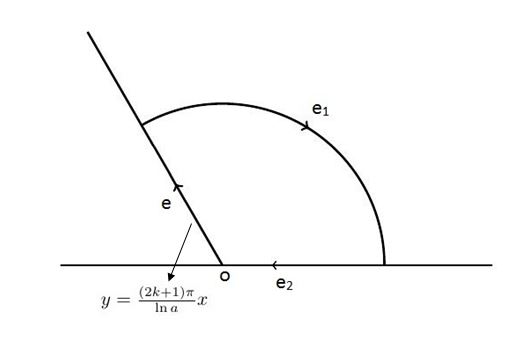}}
		\caption[caption]{Contour for the first integral}
	\end{figure}
	
	\begin{figure}[htb!]
		\centering
		\includegraphics[scale=0.6]{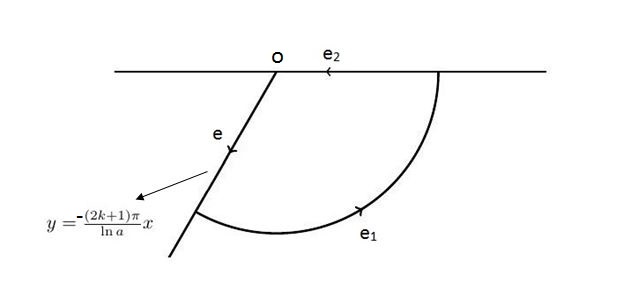}
		\caption{Contour for the second integral}
	\end{figure}

	In this case, $I_3(s) = \frac{-2s}{\pi} \sum_{k=0}^{\infty} \frac{\pi^s}{(2k+1)^{1-s} }I_{5k}(s) = \frac{2s}{\pi} \sum_{k=0}^{\infty} \frac{\pi^s}{(2k+1)^{1-s}((\frac{\ln(a)}{(2k+1)\pi})^2 + 1) }\Gamma(-s)$.
	
	$I_1(s) = \frac{\ln(a)}{s-1}I_3(s-1)  = \frac{2\ln(a)}{\pi} \sum_{k=0}^{\infty} \frac{\pi^{s-1}}{(2k+1)^{2-s}((\frac{\ln(a)}{(2k+1)\pi})^2 + 1) }\Gamma(-s + 1)$

	Thus 
	\begin{eqnarray*}
		\gamma(s) &= &  I_1 (s) + I_2(s) + I_3(s) + I_4(s) \\
		&=&  \frac{2\ln(a)}{\pi} \sum_{k=0}^{\infty} \frac{\pi^{s-1}}{(2k+1)^{2-s}((\frac{\ln(a)}{(2k+1)\pi})^2 + 1) }\Gamma(-s + 1)\\
		&+&     -\frac{(-\ln(a))^s }{2}~\Gamma (-(s-1),-\ln a)  \\
		&+&  \frac{2s}{\pi} \sum_{k=0}^{\infty} \frac{\pi^s}{(2k+1)^{1-s}((\frac{\ln(a)}{(2k+1)\pi})^2 + 1) }\Gamma(-s)\\
		&+&\frac{s}{2} (-\ln a)^s~\Gamma(-s,-\ln a)
	\end{eqnarray*}  
	
	Let $N_k = (2k+1)$ and $A = \frac{-ln(a)}{\pi}$. Let us restrict ourselves to  $e^{-\pi} < a <1$ such that 0< $\frac{A}{N_k} < 1$ for all $k$. Let $\tilde{\zeta} (s) = \sum_{k=0}^{\infty} \frac{1}{(2k+1)^s}$.
	
	Notice that 
	\begin{eqnarray*}
		\sum_{k=0}^{\infty} \frac{1}{(2k+1)^{2-s}((\frac{\ln(a)}{(2k+1)\pi})^2 + 1) } 
		&=& \sum \frac{1}{N_k^{2-s} (\frac{A^2}{N_k^2 }+ 1) }\\
		&=& \sum_{k=0}^{\infty} \frac{1}{N_k^{2-s} }[ 1- \frac{A^2}{N_k^2} + \frac{A^4}{N_k^4} - \frac{A^6}{N_k^6} ...]\\
		&=&  \tilde{\zeta}(2-s) - A^2 \tilde{\zeta}(4-s) +... \\
		&=& \sum_{n=0}^{\infty} \tilde{\zeta}((2n+2 -s) (-1)^n  A^{2n}.
	\end{eqnarray*}
	and
	\begin{eqnarray*}
		\sum_{k=0}^{\infty} \frac{1}{(2k+1)^{1-s}((\frac{\ln(a)}{(2k+1)\pi})^2 + 1) }
		&=& \sum \frac{1}{N_k^{1-s} (\frac{A^2}{N_k^2 }+ 1) }\\
		&=& \sum_{k=0}^{\infty} \frac{1}{N_k^{1-s} }[ 1- \frac{A^2}{N_k^2} + \frac{A^6}{N_k^6} - \frac{A^8}{N_k^8} ...]\\
		&=&  \tilde{\zeta}(1-s) - A^2 \tilde{\zeta}(3-s) +... \\
		&=& \sum_{n=0}^{\infty} \tilde{\zeta}(2n+1 -s) (-1)^n A^{2n}.
	\end{eqnarray*} 
	
	Thus, the following proposition gives us a functional relationship between the Dirichlet series and $\tilde{\zeta}$ function.
	
	\begin{proposition}
		Let $A = \frac{-ln(a)}{\pi}$ and $\tilde{\zeta} (s) = \sum_{k=0}^{\infty} \frac{1}{(2k+1)^s}$.
		For $\rm{Re}(s) <0$ and $e^{-\pi} < a < 1$, we have the following functional equation.
		\begin{eqnarray*}
			\gamma(s)
			&=&  \frac{2\ln(a) \pi^{s-1} \Gamma (-s+1)}{\pi}  \sum_{n=0}^{\infty} \tilde{\zeta}((2n+2 -s) (-1)^n  A^{2n}\\
			&-&     \frac{(-\ln(a))^s }{2}~\Gamma (-(s-1),-\ln a)  \\
			&+&  \frac{2s \pi^s}{\pi} \Gamma(-s) \sum_{n=0}^{\infty} \tilde{\zeta}((2n+1 -s) (-1)^n  A^{2n}\\
			&+&\frac{s}{2} (-\ln a)^s \Gamma(-s,-\ln a)
		\end{eqnarray*}
	\end{proposition}
	The infinite sums converge since $0<A<1$.
	
	\subsection{Essential singularity at $\infty$}
	
	We consider the limit of our function as $s\rightarrow \infty$ along two directions. $L_k(s,a)=\sum_{n=1}^{\infty} (-1)^{n-1}\frac{a^n}{n^s},~where~a=\frac{x^2}{c_k^2 -x^2},0<a<1,~s=\sigma+it \in \mathbb{C}$.
	\\
	\\
	\\
	When $\rm{Im}(s)$ is fixed and $\rm{Re}(s)\rightarrow \infty$ : 
	\begin{align*}
	|\sum_{n=2}^{\infty} (-1)^{n-1}\frac{a^n}{n^s}|
	&\leq \sum_{n=2}^{\infty} \frac{1}{n^{\sigma}}=\sum_{n=2}^{\infty} \frac{1}{n^{\sigma-c+c}}=\sum_{n=2}^{\infty} \frac{1}{n^{\sigma-c}}\frac{1}{n^{c}}\leq \frac{1}{2^{\sigma-c}}\sum_{n=2}^{\infty} \frac{1}{n^{c}} \rightarrow 0~as~\sigma \rightarrow \infty
	\end{align*}
	Therefore,~$\sum_{n=1}^{\infty} (-1)^{n-1}\frac{a^n}{n^s} \rightarrow a$ as $\sigma \rightarrow \infty$.
	\\
	\\
	\\
	When $\rm{Re}(s)$ is fixed and $\rm{Im}(s)\rightarrow \infty$ :
	\\
	\\
	If possible, let $L_k(\frac{1}{2}+it,a)\rightarrow a$ as $t\rightarrow \infty$. Then given $\epsilon > 0, \exists t_0 \in {\mathbb R}$ such that 
	\begin{align*}
	&\hspace{6mm}|\sum_{n=1}^{\infty} (-1)^{n-1}\frac{a^n}{n^{\frac{1}{2}+it}}-a| \leq \epsilon~~\forall~t\geq t_0
	\\
	&~or~||\sum_{n=1}^{\infty} (-1)^{n-1}\frac{a^n}{n^{\frac{1}{2}+it}}|-|a|| \leq \epsilon~~\forall~t\geq t_0
	\\
	&~or~|a|-\epsilon \leq |\sum_{n=1}^{\infty} (-1)^{n-1}\frac{a^n}{n^{\frac{1}{2}+it}}|~~\forall~t\geq t_0
	\\
	&~or~|a|-\epsilon \leq \zeta(\frac{1}{2})
	\\
	&~or~|a|\leq -1.46..+\epsilon
	\end{align*}
	which is not true for sufficiently small $\epsilon$.
	\\
	\\
	Thus we arrive at a contradiction.  The limit of $L_k(s,a)$ as $s\rightarrow \infty$ does not exist and the series has an essential singularity at infinity.

	\section{Acknowledgement:} The first author would like to thank Professor Sanoli Gun of IMSC, Chennai, India   for interesting discussions.


\begin{thebibliography}{1}
		
		\bibitem{example2} E.~Carraro, {\it Analytic continuation and functional equation of the Dirichlet eta function,\/}~~www.academia.edu.
		
		\bibitem{dey} R.~Dey: \textit{Ramanujan's identities, minimal surfaces and solitons}, Proc. Indian Acad. Sci. (Math. Sci.), 126, No. 3, (2016), 421--431.
		
		

      \bibitem{lopez1} R.~L\'opez and T.~Hasanis: \textit{Classification and construction of minimal translation surfaces in Euclidean space}, https://arxiv.org/pdf/1809.02759.pdf
		
		\bibitem{nitsche} J.C.C.~Nitsche: \textit{Lectures on Minimal surfaces} (English edition), Cambridge University Press (1989).
		
		\bibitem{osserman} R.~Osserman: \textit{Survey of minimal surfaces}, Dover Publications, New York (1986).
		
		\bibitem{ramanujan} S.~Ramanujan: \textit{Ramanujan's Notebooks} (edited by Bruce C. Berndt) (2nd ed.), Part I, Chapter 2.
		
		\bibitem{scherk} H.~F.~Scherk, \textit{Bemerkungen  \"uber die kleinste Fl\"ache innerhalb gegebener
			Grenzen}, J. Reine Angew. Math. 13 (1835), 185--208.
		
		
	\end{thebibliography}
\end{document}